\documentclass[12pt]{amsart}
\usepackage[normalem]{ulem}
\usepackage{lineno}
\usepackage{soul}
\usepackage{graphicx}
\usepackage{xcolor}
\usepackage{amssymb, amsmath, amsthm, hyperref, euscript, color}
\usepackage[dvipsnames]{xcolor}
\usepackage{amscd}
\usepackage{tikz-cd}
\usepackage{bbm}
\usepackage{enumitem}
\usepackage[english]{babel}
\usepackage{mathtools}
\usepackage{makecell}

\numberwithin{equation}{section}

\newtheoremstyle{theor}{6pt plus 1pt minus 1pt}{6pt plus 1pt minus 1pt}{\slshape}{}{\bfseries}{.}{5pt plus 1pt minus 1pt}{}
\newtheoremstyle{def}{6pt plus 1pt minus 1pt}{6pt plus 1pt minus 1pt}{}{}{\bfseries}{.}{5pt plus 1pt minus 1pt}{}
\newtheoremstyle{rmk}{6pt plus 1pt minus 1pt}{6pt plus 1pt minus 1pt}{}{}{\bfseries}{.}{5pt plus 1pt minus 1pt}{}
\newtheoremstyle{claim}{6pt plus 1pt minus 1pt}{6pt plus 1pt minus 1pt}{\slshape}{}{\bfseries}{.}{5pt plus 1pt minus 1pt}{}

\theoremstyle{theor}
\newtheorem{newstatement}{newstatement}[section]
\newtheorem{lemma}[newstatement]{Lemma}
\newtheorem{theorem}[newstatement]{Theorem}

\newtheorem{proposition}[newstatement]{Proposition}

\theoremstyle{def}
\newtheorem{definition}[newstatement]{Definition}

\theoremstyle{rmk}
\newtheorem{remark}[newstatement]{Remark}
\newtheorem*{remark*}{Remark}
\newtheorem{example}[newstatement]{Example}
\newtheorem*{example*}{Example}

\theoremstyle{claim}
\newtheorem*{claim}{Claim}

\theoremstyle{theor}
\newtheorem{thm}{Theorem}

\expandafter\let\expandafter\oldproof\csname\string\proof\endcsname
\let\oldendproof\endproof
\renewenvironment{proof}[1][\proofname]{%
  \oldproof[\slshape #1]%
}{\oldendproof}
\def\provedboxcontents#1{$\square$}

\makeatletter
\newsavebox\myboxA
\newsavebox\myboxB
\newlength\mylenA

\newcommand*\xoverline[2][0.75]{%
    \sbox{\myboxA}{$\m@th#2$}%
    \setbox\myboxB\null
    \ht\myboxB=\ht\myboxA%
    \dp\myboxB=\dp\myboxA%
    \wd\myboxB=#1\wd\myboxA
    \sbox\myboxB{$\m@th\overline{\copy\myboxB}$}
    \setlength\mylenA{\the\wd\myboxA}
    \addtolength\mylenA{-\the\wd\myboxB}%
    \ifdim\wd\myboxB<\wd\myboxA%
       \rlap{\hskip 0.5\mylenA\usebox\myboxB}{\usebox\myboxA}%
    \else
        \hskip -0.5\mylenA\rlap{\usebox\myboxA}{\hskip 0.5\mylenA\usebox\myboxB}%
    \fi}
\makeatother

\newcommand{\R}{\mathbb{R}}
\newcommand{\Z}{\mathbb{Z}}
\newcommand{\N}{\mathbb{N}}

\newcommand{\defeq}{\vcentcolon=}
\newcommand{\Bd}{\partial}
\DeclareMathOperator{\Int}{Int}
\DeclareMathOperator{\id}{id}
\newcommand{\cs}{\mathbin{\#}}
\newcommand{\bcons}{\mathbin{\natural}}

\newcommand{\CP}{{\mathbb C\mkern-0.5mu\mathrm P}}

\newcommand{\RP}{{\mathbb R\mkern-0.5mu\mathrm P}}

\DeclareMathOperator{\Pd}{PD}

\newcommand{\simtimes}{\mathbin{\widetilde{\smash{\times}}}}

\let\geq\geqslant
\let\leq\leqslant

\newcommand{\trans}[2]
  {({\text{\Large$\scriptstyle#1$}}\mkern8mu{\text{\Large$\scriptstyle#2$}})}

\newcommand{\darrow}[1]
{\mathchoice
{\xlongrightarrow{\hbox{\raisebox{0pt}[0pt][0pt]
  {\smash{$\scriptstyle\mkern1mu#1
  \mkern3mu:\mkern2mu1\mkern1.5mu$}}}}}
{\xlongrightarrow{\hbox{\raisebox{-1.5pt}[0pt][0pt]
  {\smash{$\scriptstyle\mkern2mu#1
  \mkern3mu:\mkern2mu1\mkern2.5mu$}}}}}
{}{}}

\tikzcdset{every label/.append style = {font=\normalsize,inner sep=0.75ex}}
\tikzcdset{every matrix/.append style = {row sep=5ex,column sep=1em}}

\begin{document}

\title[Branched coverings of non-orientable 4-manifolds]{Branched covering representation\\of non-orientable 4-manifolds}

\author{Valentina Bais}
\address{Scuola Internazionale Superiore di Studi Avanzati (SISSA), Via Bo\-nomea 265\\34136\\Trieste\\Italy.}
\email{vbais@sissa.it}

\author{Riccardo Piergallini}
\address{Scuola di Scienze e Tecnologie,
Università di Camerino\\Via Madon\-na delle Carceri 9\\62032 Camerino\\Italy.}
\email{riccardo.piergallini@unicam.it}

\author{Daniele Zuddas}
\address{Dipartimento di Matematica, Informatica e Geoscienze\\Università degli Studi di Trieste \\Via Valerio 12/1\\34127\\Trieste\\Italy.}
\email{dzuddas@units.it}



\begin{abstract}
We show that every closed connected non-orientable PL $4$-manifold $X$ is a simple branched covering of $\RP^4$. We also show that $X$ is a simple branched covering of the twisted $S^3$-bundle $S^1 \simtimes S^3$ if and only if the first Stiefel--Whitney class $w_1(X)$ admits an integral lift. In both cases, the degree of the covering can be any integer $d \geq 4$, provided that $d$ has the same parity as the Stiefel--Whitney number $w_1^4[X]$ in the case of $\RP^4$. Moreover, the branch set can be assumed to be non-singular if $d \geq 5$ and to have just nodal singularities if $d=4$. 
\end{abstract}

\maketitle

\kern-6pt

\section{Introduction}

A classical theorem of Alexander \cite{A1920} states that every closed PL $n$-manifold can be represented as a branched covering of $S^n$. However, no control is given over the degree and on the regularity of the branch set.
A similar result holds in the non-orientable setting for even dimensional manifolds. Indeed, Berstein--Edmonds showed in \cite[Appendix]{BE} that every closed connected non-orientable $2n$-dimensional PL manifold is a branched covering of $\RP^{2n}$, for every $n \in \mathbb{N}$. 

In dimension $3$ stronger representation theorems hold. More precisely, it was independently shown by Hilden \cite{Hi}, Hirsch \cite{Hirsch} and Montesinos \cite{M} that every closed connected orientable $3$-manifold is a simple $3$-fold covering of $S^3$ branched over a knot. Moreover, Berstein--Edmonds showed in \cite[Corollary 8.10]{BE} that every closed connected non-orientable $3$-manifold is a simple covering of $S^1 \times \RP^2$ of degree less than or equal to $6$ branched over a link. 

Moving to dimension $4$, Piergallini--Iori \cite{P,PI} proved that every closed connected orientable PL 4-manifold $X$ can be represented by a simple branched covering $p\colon X \to S^4$, whose degree $d$ and branch set $B_p \subset S^4$ satisfy one of the following conditions:

\begin{enumerate}[wide=0pt,leftmargin=*]
\item[(a)]  $d = 4$ and $B_p$ is a locally flat self-transversely immersed PL surface with no triple points,

\item[(b)] $d \geq 5$ and $B_p$ is a locally flat embedded PL surface.
\end{enumerate}




Furthermore, criteria for the existence of a simple branched covering $p\colon X \to Y$ satisfying (a) or (b) with $Y = \CP^2, S^2 \times S^2, S^2 \simtimes S^2$ were\break given by Piergallini--Zuddas \cite{PZ2} in terms of the Betti numbers $b_2^\pm(X)$.

\medskip

Here, we address the representation of closed connected non-orient\-able PL $4$-manifolds as branched covers of the basic non-orientable\break 4-manifolds $\RP^4$ and $S^1 \simtimes S^3$. The proof strategy consists in decomposing the involved manifolds into suitable orientable pieces. This allows us to exploit tools from the orientable context in order to construct branched coverings between them. Then, the final branched covering is assembled by gluing together these orientable coverings in a compatible way, using extension results due to Piergallini--Zuddas, see Theorem \ref{thm:pz} and Theorem \ref{prop:p}. 

\medskip

We work in the PL category. However, our results hold in the DIFF category as well, as every pair $(X,S)$ consisting of a PL $4$-manifold $X$ and a locally flat PL submanifold $S$ can be smoothed uniquely up to smooth isotopy.


\medskip

Before stating our main results, we recall the definition of branched covering and some related terminology for the convenience of the reader.

\medskip

Given two compact PL $n$-manifolds $X$ and $Y$, a PL map
    \[p\colon X \darrow{d} Y\]
    is called a $d$-fold \textsl{branched covering} if it is non-degenerate and restricts to a $d$-fold ordinary covering over the complement of a codimension two closed subcomplex of $Y$. The \textsl{branch set} of $p$ is the smallest subcomplex $B_p \subset Y$ with this property, while the \textsl{degree} $d=d(p)$ of $p$ is the cardinality of the preimage of any point in $Y \setminus B_p$.

    Moreover, $p$ is called a \textsl{simple} branched covering, if the preimage of a generic point of $B_p$ has cardinality $d-1$, with exactly one singular point of $p$, at which $p$ has local degree $2$.

    A $d$-fold branched covering as above is completely determined up to PL homeomorphisms by the pair $(Y,B_p)$ and the \textsl{monodromy} representation 
    \[\omega_p\colon \pi_1(Y \setminus B_p) \longrightarrow S_d\] 
    associated to the ordinary covering $p|\colon X \setminus p^{-1}(B_p) \darrow{d} Y \setminus B_p$, where $S_d$ denotes the permutation group of $\{1,\dots,d\}$. 

If $Y$ is simply connected, the fundamental group $\pi_1(Y \setminus B_p)$ is generated by a suitable set of meridians of $B_p$ and the monodromy representation $\omega_p$ is determined by the assignment of a permutation to each of these meridians so that the relations in $\pi_1(Y \setminus B_p)$ are satisfied. For $\dim Y \leq 4$ this is usually encoded by choosing a Wirtinger set of meridians with respect to some projection and by labeling each part of the projection diagram of $B_p$ with the monodromy of the corresponding meridian.

However, if $Y$ is not simply connected 
the situation is slightly more complicated, since meridians of $B_p$ do not generate $\pi_1(Y\setminus B_p)$. 
When $Y=\RP^n$, we suppose that $B_p$ is disjoint from the standard $\RP^1\subset \RP^n$ and fix the base point of $\pi_1(\RP^n \setminus B_p)$ therein. Moreover, we can assume without loss of generality that the intersection of $B_p$ with a closed tubular neighborhood of the standard $\RP^{n-1}\subset \RP^n$ is a union of fibers along the normal direction. Viewing $\RP^n$ as the union of $D^n$ with that tubular neighborhood, we describe the monodromy $\omega_p$ by its values on $[\RP^1]$ and a generating set of meridians for $B_p \cap D^n$. If $n\leq 4$, we choose a Wirtinger set of meridians for $B_p\cap D^n$ with respect to some projection.

Moreover, if $B_p$ is singular then the local monodromy at any singular point of $B_p$ is subject to further constraints, in order to guarantee that the covering space $X$ is a PL manifold.

In all cases, $p$ is a simple branched covering if and only if the monodromy of each meridian of $B_p$ is a transposition.

We observe that for every branched covering $p\colon X \to Y$ the equality \[w_1(X) = p^*(w_1(Y))\] holds. This follows from the naturality of the Stiefel--Whitney classes over $Y \setminus B_p$, noticing that removing a codimension two subcomplex of the manifold does not change $w_1$. Recall that $w_1$ is the only $1$-cohomolo\-gy class whose Kronecker pairing determines the orientability of a loop.

\medskip

Now we state our main theorems, which will be proved in Sections \ref{proofA/sec} and \ref{proofB/sec}, respectively.

\begin{thm}\label{ThmA}
   Given a closed connected non-orientable PL $4$-manifold $X$ and $d \geq 4$, there exists a $d$-fold branched covering $p\colon X \rightarrow \RP^4$ if and only if $d \equiv \langle w_1(X)^4,[X]\rangle \pmod{2}$. Moreover, we can assume that $p$ is simple and satisfies one of the above conditions \textup{(a)} or \textup{(b)}.  

\end{thm}

\begin{thm}\label{ThmB}
    Given a closed connected non-orientable PL $4$-manifold $X$ and $d \geq 4$, there exists a $d$-fold branched covering $p\colon X \rightarrow S^1 \simtimes S^3$ if and only if $w_1(X)$ admits an integral lift. Moreover, we can assume that $p$ is simple and satisfies one of the above conditions \textup{(a)} or \textup{(b)}.
\end{thm}



We expect that Theorems \ref{ThmA} and \ref{ThmB} can be generalized to branched coverings of the connected sums $\cs_m \RP^4 \cs_n S^1 \simtimes S^3$ with $m+n \geq 1$, in the spirit of \cite[Theorem 1.2]{PZ2}. We guess that the existence of such branched coverings can also be expressed in terms of certain algebraic constraints.


\subsection*{Acknowledgments}

We thank the anonymous referees for their suggestions, which substantially improved the manuscript.

Daniele Zuddas is a member of GNSAGA – Isti\-tuto Nazionale di Alta Matematica ‘Francesco Se\-veri’, Italy.



\section{Notation and conventions} 

All manifolds, submanifolds and maps between manifolds are assumed to be PL unless otherwise specified. In particular the symbol $\cong$ stands for PL homeomorphism.

$D^n$ will always denote an $n$-dimensional closed disk. $\RP^n$ is the $n$-dimensional real projective space and $S^1 \simtimes S^3$ is the twisted $S^3$-bundle over $S^1$.

For any locally flat submanifold $S \subset X$, we will indicate by $T(S,X)$ a closed regular neighborhood of $S$ inside $X$. In the following we will always have $\dim S = \dim X - 1$ or $\dim X \leq 4$, in which cases $T(S,X)$ is a $D^k$-bundle over $S$ with $k = \dim X - \dim S$, namely it is equivalent to the disk-bundle associated to $\nu_X S$, the normal bundle of $S$ in $X$. Hence, it makes sense to call $T(S,X)$ a closed tubular neighborhood, to talk about the projection $T(S,X) \to S$, and to say that $T(S,X)$ is trivial or non-trivial.

Given a closed $n$-manifold $X$, $w_i(X) \in H^i(X;\Z_2)$ denotes its $i^\text{th}$ Stiefel--Whitney class while $\Pd\colon H^i(X;G) \rightarrow H_{n-i}(X;G)$ stands for the Poincaré duality isomorphism with coefficient group $G$ with respect to a $G$-orientation (usually $G= \Z_2$ or $\Z$). An element of $H^i(X;\Z_2)$ is said to admit an integral lift if it lies in the image of the coefficient homomorphism $H^i(X;\Z) \rightarrow H^i(X;\Z_2)$.

\medskip

\section{Decomposing a non-orientable manifold}

Here, we recall a standard way to decompose a non-orientable manifold into suitable pieces. Such a decomposition turns out to be useful for constructing branched coverings of non-orientable manifolds using tools from the orientable setting. We state the following lemma in the smooth category. As noticed in the introduction, this result has a PL counterpart for every $n \leq 4$.

\begin{lemma}\label{decomposition}
    For any closed connected non-orientable smooth $n$-mani\-fold $X$, there is a closed connected orientable smooth $(n-1)$-submani\-fold $Y \subset X$ such that $[Y] = \Pd(w_1(X))$ and $X' = X \setminus \Int T(Y,X)$ is also connected and orientable. Moreover, $Y$ can be chosen in such a way that $T(Y,X)$ is trivial if and only if $w_1(X) \in H^1(X;\Z_2)$ admits an integral lift.
\end{lemma}

\begin{proof}
    Let $Y \subset X$ be a closed $(n-1)$-manifold representing the $\Z_2$-Poincaré dual of the first Stiefel--Whitney class of $X$ (cf. \cite[Théorème II.26]{T}), i.e.
    \[[Y] = \Pd(w_1(X))\in H_{n-1}(X;\Z_2).\]
    Up to a standard tubing argument, we can assume that $Y$ is connected.\break
     Now, since $[Y]=\Pd(w_1(X))$ and $X$ is non-orientable, it follows that $X'= X \setminus \Int T(Y,X)$ is a connected orientable $n$-manifold.
    
    Next, we prove that $Y$ is orientable and we consider two cases depending on the triviality of the tubular neighborhood $T(Y,X)$.
    
    If $T(Y,X)$ is trivial, then $\Bd X' = \Bd T(Y,X)$ has two boundary components that are diffeomorphic to $Y$. This implies that $Y$ is orientable.
    
   If instead $T(Y,X)$ is non-trivial, then $\Bd T(Y,X)$ is connected and the projection $\Bd T(Y,X) \to Y$ is a non-trivial (unbranched) double covering. Therefore, its non-trivial deck transformation is an involution $\iota$ of $\Bd T(Y,X) = \Bd X'$ and we have $X \cong X'/\iota$. Since $X'$ is orientable and $X$ is not, $\iota$ is necessarily orientation-preserving. This implies that $Y \cong \Bd X'/\iota$ is also orientable.

 We are left to show that $Y$ can be taken with trivial tubular neighborhood if and only if $w_1(X)$ has an integral lift. 
    
    Recall that $S^1$ is a $K(\Z,1)$ and so there is a canonical bijection $H^1(X;\Z)\cong [X,S^1]$ induced by the pullback. This implies that, whenever there is a class $\widetilde w_1(X) \in H^1(X;\Z)$ such that $\widetilde w_1(X) \equiv w_1(X) \pmod{2}$, one can find a smooth map $f\colon X \rightarrow S^1$ such that $f^*(\sigma) = \widetilde w_1(X)$, where $\sigma \in H^1(S^1;\Z) \cong \Z$ is a preferred generator. The preimage of any regular value $v \in S^1$ will give an $(n-1)$-manifold $Y \subset X$ with trivial tubular neighborhood, which we can assume to be connected up to tubing. Moreover, $Y$ is orientable since it represents $\Pd(w_1(X)) \in H_{n-1}(X;\Z_2)$. 
    
    For the converse, suppose that $Y$ has a trivial tubular neighborhood.
    An integral lift of $w_1(X)$ is represented by the map $f\colon X \rightarrow S^1$ obtained from the projection $\pi\colon T(Y,X) \cong Y \times D^1 \to D^1$ by composing with the quotient map $D^1 \to D^1/\{\pm 1\} \cong S^1$ and extending this map by the constant value $[\pm 1]$ on $X\setminus T(Y,X)$.
\end{proof}

\begin{remark}
    In Lemma \ref{decomposition}, $\Bd X' = \Bd T(Y,X)$ has two or one components depending on whether $T(Y,X)$ is trivial or not, respectively. In both cases, there exists an orientation-preserving involution $\varphi$ on $\Bd X'$ such that $X \cong X'/{\varphi}$. The map $\varphi$ is the non-trivial deck transformation of the $2$-fold covering $\Bd X' = \Bd T(Y,X) \rightarrow Y$ given by the restriction of the projection map $T(Y,X) \rightarrow Y$. In particular, if $T(Y,X)$ is trivial then $\Bd X'$ is oriented as the boundary of $X'$ with any orientation and  $\varphi$ interchanges its two components.
\end{remark}

In the proof of Theorem \ref{ThmA} (resp. \ref{ThmB}), we will decompose the $4$-mani\-fold $X$ into pieces as in Lemma \ref{decomposition}. After that, the branched covering $p\colon X \rightarrow \RP^4$ (resp. $S^1 \simtimes S^3$) will be constructed starting with a branched covering $Y \rightarrow \RP^3$ (resp. $S^3$) and extending it first to $T(Y,X)$ and then to the whole ambient $4$-manifold $X$. 

The following Lemma is a technical tool for the proof of Theorem \ref{ThmA}.

\begin{lemma}\label{lem:tubing}
    Let $X$ be a PL manifold of dimension $n \geq 3$ and let $Y\subset X$ be a closed con\-nected orientable locally flat PL $(n-1)$-submanifold with trivial normal bundle. Suppose that there exists a simple loop $\gamma\subset X$ intersecting $Y$ transversely at a single point. Then there exists a closed connected locally flat PL $(n-1)$-submanifold $Y' \subset X$ such that $[Y']=[Y]$ in $ H_{n-1}(X;\Z_2)$ and $T(Y',X)$ is not trivial. Moreover, $Y'$ can be chosen to be orientable (resp. non-orientable) if $\gamma$ is orientation-reversing (resp. orientation-preserving).
\end{lemma}
    
\begin{figure}[b]
\includegraphics{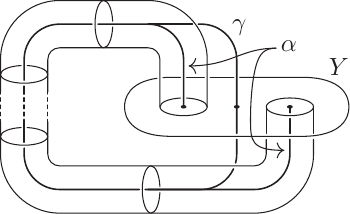}
\caption{Tubing $Y$ along the arc $\alpha$ to obtain $Y'$.}\label{tube}
\end{figure}

\begin{proof}
We think of $\gamma$ as a path with end points in $Y$ and we perturb it while keeping the end points in $Y$. In this way, we get a simple arc $\alpha \subset X$ joining two distinct points of $Y$ and not intersecting $Y$ elsewhere. Tubing $Y$ along $\alpha$ gives a closed connected PL $(n-1)$-manifold $Y'\subset X$ homologous to $Y$ in $H_{n-1}(X;\Z_2)$ (see Figure \ref{tube}). Notice that $Y'$ is PL homeomorphic to the connected sum $Y \cs (S^{n-2}\times S^1)$ (resp. $Y \cs (S^{n-2}\simtimes S^1)$) if $\gamma$ is orientation-reversing (resp. orientation-preserving). The claim on the orientability of $Y'$ follows.

We now argue on the non-triviality of $\nu_X Y'$. Let $\gamma'\subset Y'$ be a simple loop corresponding to $\{*\}\times S^1$ in $Y \cs (S^{n-2}\times S^1)$ (resp. $Y \cs (S^{n-2} \simtimes S^1)$). The conclusion follows by noticing that the restriction of $\nu_X Y'$ to $\gamma'$ is the rank $1$ twisted bundle.
\end{proof}

\section{Stabilizing branched coverings over \texorpdfstring{$\RP^n$}{RPn}} Recall that a $d$-fold branched covering $p\colon X \rightarrow S^n \text{ (resp. $D^n$})$ can be stabilized to a $(d+1)$-fold branched covering $p'\colon X \rightarrow S^n \text{ (resp. $D^n$})$ by adding an extra trivial sheet. In terms of the labeled branch set, this means that $B_{p'}$ is obtained from $B_p$ by adding a separated unknotted $(n-2)$-sphere (resp. proper trivial $(n-2)$-disk) with monodromy $\trans{i}{d+1}$ for some $i \in \{1,\dots,d\}$. The situation is rather different in the case of branched coverings over $\RP^n$, as explained here below.

\medskip

In the following, for $n \geq 2$ we will identify $T(\RP^{n-1},\RP^n)$ with $\RP^n \setminus \Int D^n$, where $D^n \subset \RP^n$ is an $n$-disk. Moreover, we note that $T(\RP^{n-1},\RP^n)$ is homeomorphic to a closed tubular neighborhood of the zero section in the non-trivial line bundle over $\RP^{n-1}$.

As already mentioned in the introduction, in this section we take base points on the standard $\RP^1$ (which we assume to be disjoint from the branch set) and describe the monodromy by its value on $[\RP^1]$ and a system of meridians for the intersection of the branch set with $D^n$.

\begin{lemma}\label{stabilization}
    Let $X$ be a compact connected $n$-manifold and
    \[p\colon X \darrow{d} T(\RP^{n-1},\RP^n)\]
    be a $d$-fold branched covering. Then there exists a $(d+2)$-fold branched covering 
    \[p'\colon X \setminus \Int D^n \darrow{d+2} T(\RP^{n-1},\RP^n),\]
    where $D^n \subset \Int X$ is a $n$-disk, the branch set $B_{p'}$ is the union of $B_p$ and a separated proper trivial $(n-2)$-disk with monodromy $\trans{i}{d+1}$ for some $i \in \{1, \dots, d\}$, and the monodromy of $\RP^1$ is given by
    \[\omega_{p'}([\RP^1]) = \omega_p([\RP^1]) \cdot \trans{d+1}{d+2}.\]
\end{lemma}

\begin{proof}
    Let $u\colon S^n \rightarrow \RP^n$ be the universal covering map and let
    \[p_1 = p \sqcup u|\colon X \sqcup T(S^{n-1},S^n) \longrightarrow T(\RP^{n-1},\RP^n)\]
    be the $(d+2)$-fold branched covering obtained by adding to $p$ two sheets with total space a $u$-equivariant tubular neighborhood $T(S^{n-1},S^n) \cong S^{n-1} \times D^1$. In particular,  the sheets of $p$ correspond to the first $d$ sheets of the cover $p_1$, while the two sheets of $u$ occupy the $(d + 1)^\text{th}$ and $(d + 2)^\text{th}$ sheets. Indicate by
    \[p_2\colon \sqcup_{d+1} D^n \longrightarrow D^n\]
    the $(d+2)$-fold branched covering of $D^n$ whose branch set is a proper trivial $(n-2)$-disk with monodromy $\trans{i}{d+1}$ for a fixed $i \in \{1, \dots, d\}$. Then, the branched covering $p'$ is defined by taking the equivariant boundary connected sum of $p_1$ and $p_2$. Namely, this operation corresponds to performing a boundary connected sum of $T(\RP^{n-1},\RP^n)$ and $D^n$ by attaching a suitable $1$-handle between them, which at the level of the covering space amounts to attaching a lift of this band for every sheet. 
    The conclusion follows by noticing that $X \bcons (S^{n-1} \times D^1) \cong X \setminus \Int D^n$.
\end{proof}

As a consequence, we get the next proposition.

\begin{proposition}\label{cor}
    Let $X$ be a closed connected $n$-manifold and 
    \[p\colon X \darrow{d} \RP^{n}\]
    a $d$-fold branched covering. Then, $p$ can be stabilized to a $(d+2)$-fold branched covering 
    \[p'\colon X \darrow{d+2} \RP^{n},\]
    where the branch set $B_{p'}$ is the union of $B_p$ and a separated trivial $(n-2)$-sphere with monodromy $\trans{i}{d+1}$ for some $i \in \{1, \dots, d\}$, and the monodromy of $\RP^1$ is given by
    \[\omega_{p'}([\RP^1]) = \omega_p([\RP^1]) \cdot \trans{d+1}{d+2}.\]
\end{proposition}

\begin{proof}
    Let $D^n \subset \RP^n$ be an $n$-disk disjoint from $B_p$, so that $p^{-1}(D^n) = \sqcup_{d} D^n$, the disjoint union of $d$ $n$-disks, and $p$ is trivial over $D^n$. Then, we consider the restriction
    \[p|\colon X \setminus \sqcup_{d} D^n \darrow{d} T(\RP^{n-1},\RP^n)\]
    and stabilize it according to Lemma \ref{stabilization} to a $(d+2)$-fold simple branched covering
    \[p_1\colon X \setminus \sqcup_{d+1} D^n \darrow{d+2} T(\RP^{n-1},\RP^n).\]
    Furthermore, we consider the $(d+2)$-fold branched covering 
    \[p_2\colon \sqcup_{d+1} D^n \darrow{d+2} D^n\]
    whose branch set is a proper trivial $(n-2)$-disk with monodromy $\trans{i}{d+1}$ for some $i \in \{1, \dots, d\}$. Notice that the choices of $i$ in the definition of $p_2$ and in the stabilization $p_1$ are the same.
    Finally, the statement follows by gluing $p_1$ and $p_2$ along the boundary, that is by putting $p' = p_1 \cup_\Bd p_2$.
\end{proof}

\begin{remark}\label{rmk:stab}
    In Proposition \ref{cor}, we have 
    \[p^* = (p')^*\colon H^*(\RP^n;\Z_2) \to H^*(X;\Z_2).\] 
    This follows from the construction of the stabilization in Lemma \ref{stabilization} and Proposition \ref{cor}, using also the fact that $d(p)\equiv d(p') \pmod{2}$.
\end{remark}

Even though we will not use this in the present paper, we notice that the proof of Proposition \ref{cor} works in a more general context. Given a $d$-fold branched covering
\[p\colon X \darrow{d} Y\]
such that $Y$ admits a $k$-fold branched covering
\[q\colon S^n \darrow{k} Y,\]
we can stabilize $p$ to a $(d+k)$-fold branched covering
\[p'\colon X \darrow{d+k} Y.\]
Namely, the branch set of $p'$ is $B_{p'} = B_p \cup B_q \cup S^{n-2}$ where $S^{n-2} \subset Y$ is a separated trivial $(n-2)$-sphere. The monodromy representation of $p'$ can be described as follows. Consider a system of generators for $\pi_1(Y\setminus (B_p \cup B_q \cup S^{n-2}))$ given by $\{[\gamma_1], \dots, [\gamma_n]\}\subset \pi_1(Y\setminus(B_p \cup B_q))$ and the class of a meridian $\mu$ of $S^{n-2}$. Then
\[\omega_{p'}([\mu])=\trans{d}{d+1} \quad \text{and} \quad \omega_{p'}([\gamma_i]) = \omega_p([\gamma_i]) \cdot \omega_q([\gamma_i])\]
for any $[\gamma_i]\in \pi_1(Y \setminus (B_p\cup B_q))$, where $\omega_q$ takes values in the permutation group of $\{d+1, \dots, d+k \}$. 

\section{Constructing branched coverings over \texorpdfstring{$T(\RP^3,\RP^4)$}{T(RP3,RP4)}}

Let us start with some basic constructions of simple branched coverings over the real projective plane.

\begin{lemma}\label{lem:coveringsRP2}
    Let $\Sigma$ be a closed connected surface and assume that $d \equiv \chi(\Sigma) \pmod{2}$. Then there is a $d$-fold simple branched covering
    \[p\colon \Sigma \darrow{d} \RP^2\]
    if either $\Sigma$ is non-orientable and $d \geq 2$ or $\Sigma$ is orientable and $d \geq 4$.
     
\end{lemma}

\begin{proof}


    If $\Sigma$ is non-orientable, we apply Lemma \ref{decomposition} to $X = \Sigma$ and set $\gamma = Y$, so that $[\gamma] = \Pd(w_1(\Sigma))$. We have a decomposition
    \[\Sigma = T(\gamma,\Sigma) \cup \Sigma'\]
    where $\Sigma' = \Sigma \setminus \Int T(\gamma, \Sigma)$ is a compact connected orientable surface. Notice that $T(\gamma,\Sigma)$ is an annulus if $\chi(\Sigma)$ is even and it is a Möbius strip otherwise.
    Since $T(\RP^1,\RP^2) \subset \RP^2$ is a Möbius strip, there is a $d$-fold cyclic unbranched covering 
    \[p_1\colon T(\gamma,\Sigma) \darrow{d} T(\RP^1,\RP^2)\]
    for every integer $d \geq 2$ with the same parity as $\chi(\Sigma)$.
    
    We want to extend this covering to a $d$-fold simple branched covering $\Sigma \rightarrow \RP^2$. In order to do so, we first observe that the orientable surface $\Sigma'$ has one or two boundary components depending on whether $d$ is odd or even. 
    In both cases, a 2-fold simple branched covering $\Sigma' \rightarrow D^2$ can be obtained as a suitable restriction of a hyper-elliptic covering of $S^2$.
    By stabilizing such a covering $d-2$ times, we get a $d$-fold simple branched covering
    \[p_2\colon \Sigma' \darrow{d} D^2.\]
    In particular, we perform the stabilizations so that the monodromy $\omega_{p_2}([\Bd D^2])$ is a cycle of length $d$ if $d$ is odd, while it is a product of two disjoint cycles of length $d/2$ if $d$ is even.
    Then, the $d$-fold simple branched coverings $p_1$ and $p_2$ can be glued together along the boundary to get a $d$-fold simple branched covering 
    \[p = p_1 \cup_\Bd p_2\colon S=T(\gamma,\Sigma)\cup \Sigma' \darrow{d} \RP^2.\]
    The gluing compatibility can be checked by observing that the monodromy $\omega_{p_1}([\RP^1])$ is a cycle of length $d$, whose square is again a cycle of length $d$ if $d$ is odd, while it is a product of two disjoint cycles of length $d/2$ if $d$ is even. To conclude, we just need to check that the domain $S$ of $p$ is $\Sigma$. In particular, it is enough to prove that $S$ is non-orientable. Indeed, since $p^{-1}(\RP^1)=\gamma$, it follows that $[\gamma]=\Pd(w_1(S))$. The conclusion follows from noticing that $\gamma$ is not trivial in homology, since the connectedness of $\Sigma'$ implies the existence of a simple loop in $S$ intersecting $\gamma$ geometrically once.

    If $\Sigma$ is orientable, one can construct a simple $d$-fold branched covering $\Sigma \darrow{d} \RP^2$ for every even $d \geq 4$ by composing a $d/2$-fold branched covering $\Sigma \darrow{d/2} S^2$ whose branch points are not antipodal with the universal covering $S^2 \darrow{2} \RP^2$. 
\end{proof}

The next lemma will be our main tool for constructing $4$-dimensional branched coverings of the tubular neighborhood of $\RP^3$ inside $\RP^4$ as fiberwise extensions of $3$-dimensional branched coverings.

\begin{lemma}\label{projcover}
    Given a closed connected orientable $3$-manifold $Y$, let $c \in H^1(Y;\Z_2)$ be any cohomology class. Then, for every $d \geq 3$ such that $d\equiv \langle c^3,[Y]\rangle \pmod{2}$, there is a $d$-fold simple covering 
    \[q\colon Y \darrow{d} \RP^3\] 
    branched over a link, such that $c=q^*(c_0)$, where $c_0$ is the non-zero element of $H^1(\RP^3;\Z_2)$.
\end{lemma}

\begin{proof}
    Let $\Sigma \subset Y$ be a closed connected surface representing the Poincar\'e dual $\Pd(c)\in H_2(Y;\Z_2)$. For every $d' \geq 3$ such that $d' \equiv \chi(\Sigma) \pmod{2}$, we will now build a $d'$-fold simple branched covering $q\colon Y \rightarrow \RP^3$ such that $q^*(c_0)=c$.  
    Consider the decompositions
    \[Y = T(\Sigma,Y) \cup Y' \quad \text{and} \quad  \RP^3 = T(\RP^2,\RP^3) \cup D^3 \]
    where $Y'= Y \setminus \Int T(\Sigma,Y)$. Lemma \ref{lem:coveringsRP2} implies the existence of a $d'$-fold simple branched covering $p\colon\Sigma \rightarrow \RP^2$. Since the ambient $3$-manifold $Y$ is orientable, $w_1(\Sigma) = w_1(\nu_Y \Sigma)$. Then, the pullback under $p$ of $\nu_0$, the normal $D^1$-bundle of $\RP^2$ in $\RP^3$, is $\nu$, the normal $D^1$-bundle of $\Sigma$ in $Y$, that is we have the following commutative diagram
    \begin{center}
    \hspace{6pt}
    \begin{tikzcd}
    T(\Sigma,Y) \arrow[d,swap,"\nu"] \arrow[rr, "q_1"] && {T(\RP^2, \RP^3)} \arrow[d, "\nu_0"] \\
    \Sigma \arrow[rr,swap, "p"] && \RP^2.
    \end{tikzcd}
    \end{center}
    By a similar logic, we also have that $w_1(\RP^2)=w_1(\nu_{\RP^3}(\RP^2))$, and hence $w_1(\nu_Y(\Sigma))=w_1(\Sigma)=p^*(w_1(\RP^2))=p^*(w_1(\nu_{\RP^3}(\RP^2))$.
    
    Here the lifting $q_1$ of $p$ is a $d'$-fold simple branched covering which fiberwise extends $p$ and $B_{q_1}=\nu_0^{-1}(B_p)$. Its restriction to the boundary is a branched covering
    \[q_1|_{\Bd}\colon \Bd T(\Sigma,Y) \darrow{d'} \Bd T(\RP^2,\RP^3) = S^2.\]
    Since the degree of $q_1|_{\Bd}$ is $d'\geq 3$, by \cite[Corollary 6.3]{BE} $q_1|_{\Bd}$ extends to a $d'$-fold simple covering
    \[q_2\colon Y' \darrow{d'} D^3\]
    branched over a proper $1$-dimensional submanifold of $D^3$. Then, the conclusion follows by setting $q = q_1 \cup_{\Bd} q_2$. Notice that $q^*(c_0)=c$ by construction and $B_q = B_{q_1}\cup B_{q_2}$ is a link in $\RP^3$. Moreover, the parity of $d'$ must coincide with that of $d \equiv \langle c^3, [Y] \rangle \pmod{2}$, since 
    \[d\equiv \langle c^3, [Y] \rangle = \langle q^*(c_0)^3,[Y] \rangle = d' \langle c_0^3,[\RP^3]\rangle \equiv d' \pmod{2}.\]
    Hence, we can take $d'=d$ in the above construction.
\end{proof}

As a straightforward consequence of Lemma \ref{projcover}, we get the following.

\begin{lemma}\label{coveringnhoodrp}
    Let $\xi\colon T \rightarrow Y$ be a $D^1$-bundle over a closed connected orientable $3$-manifold $Y$. For every $d \geq 3$ such that $d \equiv \langle w_1(\xi)^3,[Y] \rangle \pmod{2}$, there is a commutative diagram
    \begin{center}
    \hspace{24pt}
    \begin{tikzcd}
    T \arrow[d,swap,"\xi"] \arrow[rr, "\widehat q"] && {T(\RP^3, \RP^4)} \arrow[d, "\xi_0"] \\
    Y \arrow[rr,swap, "q"] && \RP^3
    \end{tikzcd}
    \end{center}
    where $q$ and $\widehat q$ are $d$-fold simple branched coverings and $\xi_0$ is the projection map of the tubular neighborhood of $\RP^3 \subset \RP^4$. In particular, $B_q\subset \RP^3$ is a link and $B_{\widehat q}=\xi_0^{-1}(B_q)$ is a non-singular proper surface.
\end{lemma}

\begin{proof}
    Apply Lemma \ref{projcover} to get a $d$-fold simple branched covering\break $q\colon Y \rightarrow \RP^3$ such that $q^*(c_0)=w_1(\xi)$. Then, the pullback of $\xi_0$ via $q$ is a $D^1$-bundle over $Y$ isomorphic to $\xi$. We can hence define $\widehat q\colon T \rightarrow T(\RP^3,\RP^4)$ to be the fiberwise extension of $q$ (or, equivalently, its pullback under $\xi_0$), which satisfies the required conditions. 
\end{proof}

\section{Fillability properties of branched coverings of \texorpdfstring{$S^3$}{S3}}

The following notion of fillability will be crucial for the construction of the branched coverings in the proof of Theorems \ref{ThmA} and \ref{ThmB}.

\begin{definition}
    A $d$-fold simple branched covering $p\colon Y \rightarrow S^3$ is said to be \textsl{ribbon fillable} if it can be extended to a $d$-fold simple branched covering $q \colon X \rightarrow D^4$ whose branch set $B_q \subset D^4$ is a properly embedded ribbon surface. This immediately implies that $\Bd X = Y$ and $\Bd B_q = B_p$ as labeled links.
\end{definition}

\begin{remark}\label{rmk:rf}
    As observed in \cite[Section 1]{PZ1}, the above definition is invariant under equivalence of $p$ up to homeomorphism and this implies that the ribbon fillability of $p$ can be expressed in terms of the labeled branch set $B_p$ by requiring that it is a labeled link in $S^3$ bounding a labeled ribbon surface in $D^4$.
\end{remark}

In the following, we will denote by $\Gamma_d$ the cobordism group of $d$-fold simple coverings of $S^3$ branched over a link. More precisely, we set
\[\Gamma_d = \{p\colon Y \to S^3 \text{ a $d$-fold simple covering branched over a link}\}/{\sim}\,,\]
where $Y$ is any closed oriented 3-manifold and $p_0\sim p_1$ if and only if there is a $d$-fold simple covering $q\colon X \rightarrow S^3 \times I$ branched over a proper non-singular surface in $S^3 \times I$, which restricts to $p_i$ over $S^3 \times \{i\}$ for $i=0,1$.
A group operation can be defined on $\Gamma_d$ by setting $[p] + [p'] = [p'']$, where $p''$ is the branched covering described by the separated union of the labeled branch links of $p$ and $p'$.

The group $\Gamma^{\textup{or}}_d$ is defined analogously, by additionally requiring that the branch set of the covering $p$ is an oriented link and the branch set of the cobordism $q$ is an oriented surface.

It is worth noting that a $d$-fold simple covering $p \colon Y \darrow{d} S^3$ branch\-ed over a link represents the trivial class in $\Gamma_d$ if and only if $B_p$ bounds a proper non-singular labeled surface $F\subset D^4$ (required to be orientable in the case of $\Gamma_d^\textup{or}$), which represents a $d$-fold simple branched covering $\widetilde p \colon X\to D^4$ such that $\Bd X = Y$, $\widetilde p\mkern2mu |_Y = p$ and $\Bd F = B_p$ as labeled links.

The relevance of ribbon fillability in our context is due to the following two results in \cite{PZ1}. In particular, Theorem \ref{thm:pz} is a rephrasing of \cite[Theorem 1.2]{PZ1}, which is one of the main results therein, while Theorem \ref{prop:p} is stated without proof as \cite[Theorem 1.8]{PZ1}. Since we will make use of the latter theorem in the present work, we provide a proof of it. We emphasise that, while in Theorem \ref{thm:pz} the $4$-manifold $X$ is required to have exactly $n$ boundary connected components, in Theorem \ref{prop:p} we allow this number to possibly be greater than $n$.

\begin{theorem}[Piergallini--Zuddas]\label{thm:pz}
    Let $X$ be a compact connected oriented PL $4$-manifold with $n\geq 1$ boundary components and let\break $p\colon \Bd X \rightarrow \sqcup_n S^3$ be a disjoint union of $d$-fold ribbon fillable simple branched coverings, with $d \geq 4$. Then, $p$ can be extended to a $d$-fold simple branched covering $q \colon X \rightarrow S^4 \setminus \sqcup_n \Int D^4$ whose branch set $B_q$ satisfies one of the above conditions \textup{(a)} or \textup{(b)}. 
\end{theorem}
 
\begin{theorem}[Piergallini--Zuddas]\label{prop:p} Let $X$ be a compact connected oriented PL $4$-manifold and let $p\colon \Bd X \rightarrow \sqcup_n S^3$ be a disjoint union of $d$-fold ribbon fillable simple branched coverings, with $d \geq 4$ and $n \geq 1$. Then, $p$ can be extended to a $d$-fold simple branched covering $q \colon X \rightarrow S^4 \setminus \sqcup_n \Int D^4$ whose branch set $B_q$ satisfies one of the above conditions \textup{(a)} or \textup{(b)}. 
\end{theorem}

\begin{proof}
Let $Y_1, \dots, Y_m$ be the connected components of $\Bd X$ for some $m \geq n$. We prove the statement by induction on $m$ starting from the case $m=n$, where the conclusion directly follows from Theorem \ref{thm:pz}. 

Suppose now that $m > n$. Denote by $C\cong \sqcup_{i=1}^m Y_i \times [0,1]$ a collar of $\Bd X$ in $X$ and by $C_i\cong Y_i \times [0,1]$ the connected component of $C$ containing $Y_i$. We number the $3$-spheres in the disjoint union as $S^3_1, \dots, S^3_n$ in such a way that over $S^3_1$ there are at least two components of $\Bd X$. Accordingly, we set
\[p_i=p|\colon p^{-1}(S^3_i)\darrow{d} S^3_i.\] 

We assume that the components $Y_i$ are numbered in such a way that\break $p^{-1}(S^3_1)=Y_1\sqcup \dots \sqcup Y_k$ for some $k \geq 2$. We thicken $p_1$ by crossing with the identity of $[0,1]$ to get a simple $d$-fold branched covering 
\[q_1 = p_1 \times \id_{[0,1]}\colon \sqcup_{i=1}^k C_i \darrow{d} S^3_1\times [0,1].\]
We number the sheets of $q_1$ so that the first and the second sheet lie in the components $C_1$ and $C_2$ respectively. Moreover, since $X$ is connected, we can find inside $X\setminus \Int(C)$ a $1$-handle $H^1$ attached to $C$ that joins $C_1$ and $C_2$. Then, we can extend $q_1$ over $H^1$ to get a simple $d$-fold branched covering
\[q_1'\colon \sqcup_{i=1}^k C_i \cup H^1 \darrow{d} S^3_1 \times [0,1]\]
having as branch set the union of $B_{q_1}$ and a separated proper trivial $2$-disk with monodromy $\trans{1}{2}$, whose boundary is contained in $S^3_1 \times \{1\}$. In other words, $q_1'$ is obtained by performing an equivariant boundary connected sum between $q_1$ and the $d$-fold branched covering\break $\sqcup_{d-1} D^4 \mapsto D^4$ whose branch set is a proper trivial $2$-disk with monodromy $\trans{1}{2}$.
By construction, the restriction of $q_1'$ over $S^3_1 \times \{1\}$ is a ribbon fillable simple $d$-fold branched covering
\[p_1' = q_1'|_{S^3_1\times \{1\}}\colon (Y_1\cs Y_2) \sqcup Y_3 \sqcup \dots \sqcup Y_k \darrow{d} S^3_1.\]  
We consider the union of $d$-fold ribbon fillable simple branched coverings
\[p' = p_1'\sqcup p_2 \sqcup \dots \sqcup p_n\colon (Y_1 \cs Y_2) \sqcup Y_3 \sqcup \dots \sqcup Y_m \darrow{d} \sqcup_{i=1}^n S^3_i\]
having as domain a $3$-manifold with $m-1$ connected components. By the inductive hypothesis, we can extend $p'$ to a simple $d$-fold branched covering
\[q_2\colon X\setminus \Int (\sqcup_{i=1}^k C_i \cup H^1) \darrow{d} S^4 \setminus \Int(\sqcup_n D^4)\]
whose branch set satisfies one of the conditions (a) or (b). The conclusion follows by setting $q = q_1' \cup_{\Bd} q_2$.
\end{proof}

The next lemma tells us that a $d$-fold simple branched covering of $S^3$ is ribbon fillable if and only if it represents the null element in $\Gamma_d$.

\begin{lemma}\label{lem:fillableimplicaribbon}
    Let $p\colon Y \rightarrow S^3$ be a connected $d$-fold simple branched covering whose labeled branch set $B_p \subset S^3$ bounds a proper non-singular surface $F \subset D^4$ labeled with transpositions in $S_d$. Then $B_p$ also bounds a ribbon surface $F'\subset D^4$ which is again labeled with transpositions in $S_d$. Moreover, $F'$ can be taken orientable if $F$ is orientable. Therefore, $p$ represents the null class in $\Gamma_d$ (resp. $\Gamma_d^\textup{or}$) if and only if $p$ is ribbon fillable (resp. ribbon fillable by an orientable ribbon surface).
\end{lemma}

\begin{proof}
    Without loss of generality, we can assume that the restriction $\rho|_{\Int F}$ of the squared norm function $\rho \defeq \| \cdot \|^2\colon D^4 \to \R$ is Morse. Then $F$ is ribbon if and only if $\rho|_{\Int F}$ has no local maxima. 

    If $F$ is not ribbon, by general position we also assume that each radius of $D^4$ containing a local maximum of $\rho|_{\Int F}$ does not meet $F$ elsewhere. Then, we proceed by subsequently eliminating the local maxima as follows. Let $k$ be the number of the local maxima of $\rho|_{\Int F}$. Pick one of them and push a small $2$-disk neighborhood of it in $F$ radially to $S^3$. In this way, we obtain a surface that intersects $S^3$ along $B_p$ and a separated $2$-disk. 
    By removing the interior of that disk we get a proper non-singular labeled surface whose boundary is the union of $B_p$ and a separated unknot labeled with a single transposition in $S_d$. That unknot can be joined to $B_p$ by a suitable labeled band in $S^3$. Indeed, since $Y$ is connected, the image of the monodromy $\omega_p$ acts transitively by conjugation on the transpositions in $S_d$. This image is generated by meridians of the branch set, then the existence of a connecting band with the right labeling follows. Finally, by pushing the interior of the band inside $D^4$ we get a proper non-singular labeled surface, which we still denote by $F$, such that $\rho|_{\Int F}$ has $k-1$ local maxima and whose boundary is label-isotopic to $B_p$. By iterating this modification over all the local maxima, we eventually obtain the desired labeled ribbon surface $F'$.

    Notice that this operation preserves the orientability of the labeled surfaces involved in the construction, although it changes its topology, since $\chi(F') = \chi(F) - 2 k$, as well as the topology of the interior of the covering 4-manifold.
\end{proof}

We recall the following result from \cite[Theorems 9 and 11]{HL}.

\begin{figure}[ht]
\includegraphics{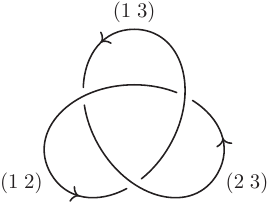}
\caption{The labeled branch set of a generator of $\Gamma_3^{\text{or}}$.}\label{trefoil}
\end{figure}

\begin{theorem}[Hilden--Little]\label{thm:cob}
    There are exactly three cobordism\break classes of $3$-fold simple coverings $p\colon Y \rightarrow S^3$ branched over oriented links, and hence $\Gamma^{\textup{or}}_3$ is cyclic of order three. A generator is given by the irregular $3$-fold simple covering $S^3 \rightarrow S^3$ branched over an oriented (left or right-handed) trefoil knot, see Figure \ref{trefoil}.
\end{theorem}

The next lemma is a direct consequence of Theorem \ref{thm:cob}.

\begin{figure}[ht]
\includegraphics{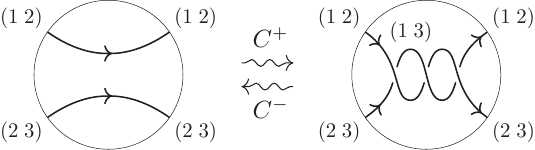}
\caption{The oriented Montesinos moves $C^{\pm}$.}\label{c}
\end{figure}

\begin{lemma}\label{lem:montesinos}
    Applying an oriented Montesinos move $C^{\pm}$ (see Figure \ref{c}) to a $3$-fold simple covering $p\colon Y \rightarrow S^3$ branched over an oriented link changes the cobordism class of the covering in $\Gamma_3^{\textup{or}} \cong \Z_3$ by $\pm 1$, while leaving $Y$ unchanged up to homeomorphism. 
\end{lemma}

\begin{proof}
That the moves $C^\pm$ do not change the manifold $Y$ is a well known fact (cf. \cite{M1, M2}). For the rest of the statement it suffices to observe that, up to oriented cobordism, performing a $C^\pm$ move is the same as adding/deleting a separated copy of the labeled knot in Figure \ref{trefoil}.
\end{proof}

\begin{remark}\label{mont}
    The oriented Montesinos moves $C^{\pm}$ can be applied, up to label conjugation in $S_d$, to any connected $d$-fold simple covering between $3$-manifolds branched over a link, with $d \geq 3$. This operation does not change 
    the homotopy class of the given branched covering, since the modifications occur locally inside $3$-disks.
\end{remark}

\section{Proof of Theorem \ref{ThmA}}\label{proofA/sec}

We start by proving the necessary condition on the degree. Given any $d$-fold branched covering $p\colon X \rightarrow \RP^4$ we have
\begin{eqnarray*}
\langle w_1(X)^4,[X] \rangle &\!\!\!=\!\!\!& \langle p^*(w_1(\RP^4))^4,[X] \rangle = \langle w_1(\RP^4)^4, p_*[X] \rangle \\ 
&\!\!\! = \!\!\!& \langle w_1(\RP^4)^4,d\,[\RP^4] \rangle \equiv d \pmod{2}.
\end{eqnarray*}
  
For the converse, let $d\geq 4$ such that $\langle w_1(X)^4,[X]\rangle \equiv d \pmod{2}$. We consider the decompositions
\[X = T(Y,X) \cup_{\Bd} X' \quad \text{and} \quad \RP^4 = T(\RP^3,\RP^4) \cup_{\Bd} D^4\]
given by Lemma \ref{decomposition}, where $Y \subset X$ is a connected orientable $3$-dimen\-sional submanifold such that 
\[[Y] = \Pd(w_1(X))\in H_3(X;\Z_2)\] and $X' = X \setminus \Int T(Y,X)$. Let $\Sigma \subset Y$ be a connected surface such that 
\[[\Sigma]=\Pd (w_1(\nu_X Y))\in H_2(Y;\Z_2).\] 
Notice that $\Sigma$ is obtained by intersecting $Y$ with a transversal copy $Y'$. Then, we have
\begin{eqnarray*}
&&\kern-20pt \langle w_1(X)^4,[X] \rangle = [Y \cap Y']^2 =[\Sigma]^2= \langle w_2(\nu_X \Sigma),[\Sigma] \rangle \\
&& = \langle w_2(\nu_Y \Sigma \oplus \nu_X Y|_{\Sigma}), [\Sigma] \rangle=\langle w_1(\nu_Y \Sigma)\smile w_1( \nu_X Y|_{\Sigma}), [\Sigma] \rangle \\
&& = \langle w_1(\nu_Y \Sigma)^2,[\Sigma] \rangle =\langle w_1(\Sigma)^2,[\Sigma]\rangle \equiv \chi(\Sigma) \pmod{2}.
\end{eqnarray*}
Here, we use the identity $w_1(\Sigma) = w_1(\nu_Y \Sigma)$, implied by the orientability of $Y$, and the isomorphism of line bundles $\nu_Y \Sigma \cong \nu_X Y|_{\Sigma}$.  The latter can be obtained by viewing $\Sigma=Y \cap Y'$ as the zero locus of a section $s \colon Y \rightarrow \nu_XY$ and by restricting its differential to $\nu_Y \Sigma$.

Also note that we have the congruence $d \equiv \langle w_1(\nu_X Y)^3, [Y]\rangle \pmod{2}$, which is implied by the fact that
\[\langle w_1(\nu_X Y)^3, [Y] \rangle = \langle w_1(X)^3,i_*([Y]) \rangle = \langle w_1(X)^4,[X] \rangle,\]
where $i_*$ is the homomorphism induced in homology by the inclusion $i\colon Y\to X$. The first equality follows from the orientability of $Y$, while the second one is a consequence of the duality between the Kronecker product and the algebraic intersection of homology classes.

We will now construct a $d$-fold simple branched covering 
\[p\colon X \darrow{d} \RP^4\]
for any $d \geq 4$ such that $d \equiv \langle w_1(X)^4,[X]\rangle \pmod{2}$.
    
We distinguish two cases, according to the parity of $d$.
    
\smallskip\noindent
\textbf{Case 1: \boldmath $d$ odd.}
By Lemma \ref{coveringnhoodrp}, there exists a $3$-fold simple branched covering
\[q\colon Y \darrow{3} \RP^3,\]
which admits a fiberwise extension to a $3$-fold simple branched covering
\[\widehat q\colon T(Y,X) \darrow{3} T(\RP^3,\RP^4).\]
Notice that the restriction of $\widehat q$ over $\Bd T(\RP^3,\RP^4)$ coincides with the pullback \[\widetilde q\colon \Bd T(Y,X) \darrow{3} \Bd T(\RP^3,\RP^4)\] of $q$ under the universal covering $u\colon S^3 \rightarrow \RP^3$, which is given by the restriction of the projection $T(\RP^3,\RP^4) \to \RP^3$ to the boundary.
    
We can suppose without loss of generality that $\widetilde q$ is ribbon fillable. Indeed, if this is not the case, Lemma \ref{lem:fillableimplicaribbon} implies that $\widetilde q$ does not represent the trivial element in the cobordism group $\Gamma_3^{\textup{or}}$ of $3$-fold simple branched coverings for any orientation of the branch set $B_{\widetilde{q}}$. Fixing an arbitrary orientation on $B_q$, we consider the $u$-invariant orientation induced on $B_{\widetilde q}$ and observe that applying one oriented Montesinos move $C^{\pm}$ to $q$ amounts to applying two oriented Montesinos moves of the same type to the pullback $\widetilde q$. Then, by Theorem \ref{thm:cob} and Lemma \ref{lem:montesinos} a single oriented Montesinos move on $q$ suffices to make $\widetilde q$ trivial in $\Gamma_3^\text{or}\cong \Z_3$, and hence ribbon fillable by Lemma \ref{lem:fillableimplicaribbon}. Notice that the domains of $\widehat q$ and $\widetilde q$ are unchanged. Moreover, Remark \ref{mont} implies that an oriented Montesinos move $C^{\pm}$ does not change the homotopy class of a given covering and the induced pullbacks are hence not affected by this operation.
    
At this point, we assume that $\widetilde q$ is ribbon fillable. We stabilize $k$ times $q$ as in Proposition \ref{cor} to get a $d$-fold simple branched covering
\[q' \colon Y  \darrow{d} \RP^3\]
with $k=(d-3)/2$. Then, we consider its fiberwise extension
\[\widehat q\,'\colon T(Y,X) \darrow{d} T(\RP^3,\RP^4),\]
whose total space is still $T(Y,X)$, since thanks to Remark \ref{rmk:stab} we have 
\[q^*(w_1(\nu_{\RP^4} \RP^3))=(q')^*(w_1(\nu_{\RP^4} \RP^3)),\] 
and the $D^1$-bundles associated to $q^*(\nu_{\RP^4} \RP^3)$ and  $(q')^*(\nu_{\RP^4} \RP^3)$ are hence isomorphic. We also consider the restriction of $\widehat q \,'$ to the boundary
\[\widetilde q\,'\colon \Bd T(Y,X) \darrow{d} \Bd T(\RP^3,\RP^4).\]  
The branch set of $q'$ is the union of $B_{q}$ and $k$ separated unknots with monodromies $\trans{3}{4}, \trans{5}{6}, \dots, \trans{d-2}{d-1}$ and thus $\widetilde q\,'$ is still ribbon fillable, given that its branch set is the union of $B_{\widetilde q}$ with $d-3$ separated unknots (in fact, it is obtained from $\widetilde q$ by stabilizing $d-3$ times). We can hence apply Theorem \ref{thm:pz}  with $n=1$ to extend $\widetilde q\,'$ and hence $\widehat q\,'$ as well to a $d$-fold simple branched covering
\[p\colon X = T(Y,X) \cup_{\Bd} X' \darrow{d} \RP^4 = T(\RP^3,\RP^4) \cup_{\Bd} D^4\]
whose branch set $B_{p}$ is a non-singular locally flat surface, satisfying condition (b).

\smallskip\noindent
\textbf{Case 2: \boldmath $d$ even.}  If $\Sigma$ is orientable, we can make it non-orientable up to homology in $H_2(Y;\Z_2)$ as follows. Lemma \ref{lem:tubing} implies that we can assume that $w_1(\nu_X Y)$ is not trivial, up to replacing $Y$ with a $3$-manifold representing the same class in $H_3(X;\Z_2)$. Since $[\Sigma]=\Pd (w_1(\nu_X Y))$ by definition, it follows that $[\Sigma] \neq 0 \in H_2(Y;\Z_2)$. An additional application of Lemma \ref{lem:tubing} implies that, up to a single non-orientable tubing, we can assume that $\Sigma$ is non-orientable.

By Lemma \ref{lem:coveringsRP2} and its proof, there is a $2$-fold simple branched covering
\[\Sigma \darrow{2} \RP^2\]
such that the standard $\RP^1\subset \RP^2$ is disjoint from the branch set and its monodromy is the transposition $\trans{1}{2}$. By the Riemann--Hurwitz formula, the number of branch points is $2-\chi(\Sigma)$, which is even by the assumption on the parity of $d$. Since $w_1(\Sigma) = w_1(\nu_Y \Sigma)$ and $w_1(\RP^2)$ pulls back to $w_1(\Sigma)$, we can extend this covering fiberwise to a $2$-fold simple branched covering
\[T(\Sigma,Y) \darrow{2} T(\RP^2, \RP^3).\]
Then, we stabilize the latter map as in Lemma \ref{stabilization} to get a $4$-fold simple branched covering
\[q_1\colon T(\Sigma,Y) \setminus \Int D^3 \darrow{4} T(\RP^2,\RP^3),\]
whose branch set consists of the preimages with respect to the projection $T(\RP^2,\RP^3) \rightarrow \RP^2$ of the $2h=2-\chi(\Sigma)$ branch points in $\RP^2$ and a proper separated trivial arc with monodromy $\trans{2}{3}$. Moreover, after such stabilization, the monodromy of the standard $\RP^1$ is given by 
\[\omega_{q_1}([\RP^1]) = \trans{1}{2}\trans{3}{4}.\] 
Notice that the restriction of $q_1$ to the boundary is a $4$-fold simple branched covering
\[q_1|_{\Bd}\colon \Bd T(\Sigma,Y) \sqcup S^2 \darrow{4} S^2\]
which is $3$-fold on $\Bd T(\Sigma,Y)$ and one to one on $S^2$. We can extend this restriction to a $4$-fold simple branched covering
\[q_2\colon (Y \setminus \Int T(\Sigma,Y))\sqcup D^3 \darrow{4} D^3\]
by applying the Berstein--Edmonds extension result \cite[Corollary 6.3]{BE} on $Y \setminus \Int T(\Sigma,Y)$ and taking the disjoint union with the identity map of $D^3$. We can then glue $q_1$ and $q_2$ along the boundary to get a $4$-fold simple branched covering
\[q=q_1 \cup_{\Bd} q_2\colon Y \darrow{4} \RP^3.\]
Since $[\Sigma] = [q^{-1}(\RP^2)]$ by construction and \[[\RP^2]=\Pd(w_1(\nu_{\RP^4}(\RP^3)),\] it follows by Poincar\'e duality that \[w_1(\nu_X Y) = q^*(w_1(\nu_{\RP^4} \RP^3)).\] Hence, we can extend $q$ fiberwise to a $4$-fold simple branched covering 
\[\widehat q\colon T(Y,X) \darrow{4} T(\RP^3,\RP^4),\]
whose restriction  
\[\widetilde q\colon \Bd T(Y,X) \darrow{4} \Bd T(\RP^3,\RP^4) = S^3\]
can be identified with the pullback of $q$ under the universal covering $u\colon S^3 \rightarrow \RP^3$. 
    
\begin{claim}
    The branched covering $\widetilde q$ can be assumed to be ribbon fillable, up to modifying $q$ by a suitable Montesinos move $C^{\pm}$. In particular, the domains of $\widehat q$ and $\widetilde q$ are unchanged, as in Case 1.
\end{claim}

\begin{figure}[b]
\includegraphics{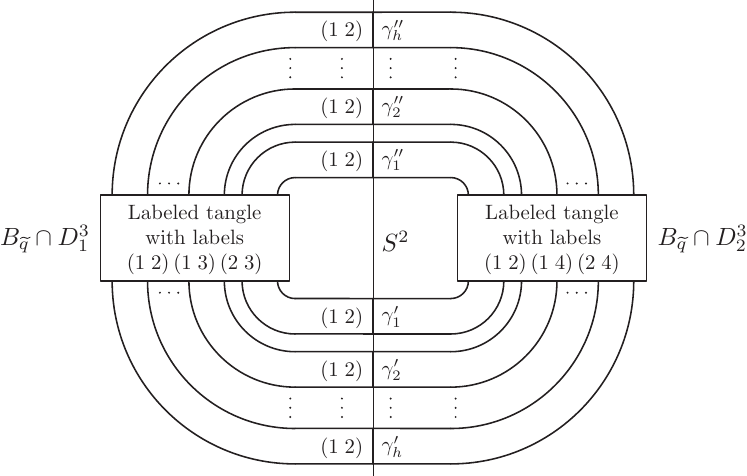}
\caption{The branch set $B_{\widetilde q} \subset \partial T(\RP^3,\RP^4) = S^3$.}\label{qtilde}
\end{figure}

\begin{proof}[Proof of the claim]
The equatorial $S^2$ in $S^3$ given by the preimage of the standard $\RP^2 \subset \RP^3$ intersects the branch set of $\widetilde q$ in exactly $4h$ points with monodromy $\trans{1}{2}$. This follows from the fact that the restriction of $\widetilde q$ to $S^2$ coincides with the pullback of the restriction of $q$ to $\RP^2$ under the universal covering $S^2 \rightarrow \RP^2$, noticing that the monodromies of the branch points of $q$ over $\RP^2 \setminus \RP^1$ are all equal to $\trans{1}{2}$. Moreover, such $S^2$ determines a splitting
    \[S^3 = D^3_1 \cup_{S^2\!} D^3_2\]
    into two 3-dimensional hemispheres. Now, since $\widetilde q$ is the pullback of $q$ under the universal covering $u \colon S^3 \rightarrow \RP^3$, its branch set $B_{\widetilde q}$ is symmetric with respect to the antipodal map $\alpha\colon S^3 \rightarrow S^3$, which acts on the set of labels by conjugation with $\trans{1}{2} \trans{3}{4}$.
    
    We now claim that the restriction of $\widetilde q$ over $D^3_1$ and $D^3_2$ has a $3$-fold connected component and a separated trivial sheet, and the labels of $B_{\widetilde q}\cap D^3_i$ take values in $\{\trans{1}{2},\trans{2}{3}, \trans{1}{3}\}$ if $i=1$ and $\{\trans{1}{2},\trans{1}{4},\break \trans{2}{4}\}$ if $i=2$ (see Figure \ref{qtilde}). Indeed, recall that the labels of $B_q$ inside $\RP^3\setminus \RP^2$ belong by construction to $\{\trans{1}{2},\trans{2}{3},\trans{1}{3}\}$.
    
    Since the universal covering $u \colon S^3 \rightarrow \RP^3$ sends $D^3_1\setminus S^2$ homeomorphically onto $\RP^3 \setminus \RP^2$ and $\widetilde q = u^*(q)$, $B_{\widetilde q}\cap D^3_1$ is also labeled by transpositions in $\{\trans{1}{2},\trans{2}{3},\trans{1}{3}\}$. In order to get $B_{\widetilde q}\cap D^3_2$ together with its labels, it is enough to apply the antipodal map to $B_{\widetilde q}\cap D^3_1$ and to transport the labels by conjugation via $\trans{1}{2} \trans{3}{4}$. In this way, $\trans{1}{2}$ is sent to itself, while $\trans{2}{3}$ and $\trans{1}{3}$ go to $\trans{1}{4}$ and $\trans{2}{4}$ respectively. 

    We also consider the decomposition 
    \[D^4 = D^4_1 \cup_{D^3} D^4_2,\] where $D^3\subset D^4$ is the standard $3$-disk spanned by $S^2\subset S^3$ and the $D^4_i$'s are two half $4$-disks, whose boundaries will be denoted by 
    \[S^3_i = \Bd D^4_i = D^3 \cup_{S^2} D^3_i\] 
    for $i=1,2$. 
    
    A labeled ribbon surface in $D^4$ bounding $B_{\widetilde q}$ can be constructed as follows (cf. Figure \ref{qtilde}). Join two by two the $4h$ branch points in $B_{\widetilde q}\cap S^2$ with an $\alpha$-invariant collection of disjoint arcs $\gamma'_1,\gamma''_1, \dots, \gamma'_h, \gamma''_h$ in $S^2$. These can be constructed by taking the preimage of disjoint arcs $\gamma_1,\dots,\gamma_h$ in $\RP^2$ connecting in pairs the (even number of) branch points of the $2$-fold cover $\Sigma \rightarrow \RP^2$. 

    For $i=1,2$, we define the labeled link 
    \[L_i = (B_{\widetilde q}\cap D^3_i) \cup \gamma'_1 \cup \gamma''_1 \cup \dots \cup \gamma'_h \cup \gamma''_h \subset S^3_i,\] 
    where the arcs $\gamma'_1,\gamma''_1 \dots \gamma'_h, \gamma''_h$ have the label $\trans{1}{2}$, so that $L_i$ represents a connected $3$-fold simple branched covering of $S^3_i$. 
    Up to performing a single oriented Montesinos move $C^{\pm}$ on $q$ as in Case 1, we can assume that $L_1$ bounds an oriented labeled ribbon surface $R_1 \subset D^4_1$, thanks to Lemma \ref{lem:montesinos} and Lemma \ref{lem:fillableimplicaribbon}. By forgetting the orientation and acting on $R_1$ with the antipodal map of $S^3$ we get a ribbon surface $R_2\subset D^3_2$ for $L_2$, whose labels are obtained from those of $R_1$ by conjugation with $\trans{1}{2} \trans{3}{4}$. Now, we have that 
    \[R_1 \cap R_2=\Bd R_1 \cap \Bd R_2= \gamma'_1 \cup \gamma''_1 \cup \dots \cup \gamma'_h \cup \gamma''_h\]
    with labels $\trans{1}{2}$. Therefore, $R \defeq R_1 \cup R_2\subset D^4$ is a labeled ribbon surface having $B_{\widetilde q}$ as the boundary.
\end{proof}

At this point, we resume the proof of Theorem \ref{ThmA} with $\widetilde q$ ribbon fillable. We can stabilize $k=(d-4)/2$ times $q $ as in Proposition \ref{cor} to get a $d$-fold simple branched covering
\[q'\colon Y \darrow{d} \RP^3,\]
whose fiberwise extension
\[\widehat q\,'\colon T(Y,X) \darrow{d} T(\RP^3,\RP^4),\]
is still defined on $T(Y,X)$ in virtue of Remark \ref{rmk:stab}, as in Case 1. Moreover, we denote by 
\[\widetilde q\,'\colon \Bd T(Y,X) \darrow{d} \Bd T(\RP^3,\RP^4)\]
the restriction of $\widehat q\,'$ to the boundary. As in the previous case, $\widetilde q\,'$ is still ribbon fillable, since its branch set is the union of $B_{\widetilde q}$ with $d-4$ separated unknots (in fact, it is obtained from $\widetilde q$ by stabilizing $d-4$ times). 
  
Theorem \ref{prop:p} with $n=1$ implies that we can extend $\widetilde q\,'$ and hence $\widehat q\,'$ as well to a $d$-fold simple branched covering
\[p\colon X=T(Y,X) \cup_{\Bd} X' \darrow{d} \RP^4=T(\RP^3,\RP^4)\cup_{\Bd} D^4\]
whose branch set $B_p\subset \RP^4$ has the desired properties (a) or (b).

\hfill \qed

\section{Proof of Theorem \ref{ThmB}}\label{proofB/sec}

We start by proving the necessary condition. Let $\widetilde w_1( S^1 \simtimes S^3)$ be any integral lift of $w_1(S^1 \simtimes S^3)$. If there is a $d$-fold branched covering $p\colon X \rightarrow S^1 \simtimes S^3$, then $p^*(\widetilde w_1(S^1 \simtimes S^3))$ is an integral lift of $w_1(X)$.

Conversely, let $d\geq 4$ be fixed and suppose that $w_1(X)$ admits an integral lift. Lemma \ref{decomposition} implies the existence of a decomposition
\[X = T(Y,X) \cup_{\Bd} X',\]
where $Y\subset X$ is a connected orientable $3$-manifold (denoted by $N$ in the lemma) having trivial tubular neighborhood and whose $\Z_2$-homol\-ogy class satisfies $[Y]=\Pd(w_1(X))$, while $X'=X \setminus \Int T(Y,X)$. As the first step in the construction of the desired branched covering $p\colon X \to S^1 \simtimes S^3$, we consider a $d$-fold ribbon fillable simple branched covering
\[q\colon Y \darrow{d} S^3.\]
This can be obtained by stabilizing $d-3$ times a $3$-fold ribbon fillable simple branched covering $Y \to S^3$, whose existence is established by Montesinos in \cite{M} or by Bobtcheva--Piergallini in \cite{BP}. Since $T(Y,X)$ is a product neighborhood, we can extend $q$ to a simple $d$-fold covering
\[\widehat q \,\cong q \times \id_{D^1}\colon T(Y,X)\cong Y \times D^1 \darrow{d} T(S^3,S^1 \simtimes S^3)\cong S^3 \times D^1\]
branched over a collection of proper annuli, where $S^3$ is a fiber of the twisted bundle $S^1 \simtimes S^3$. This restricts to the boundary as a $d$-fold simple branched covering
\[\widetilde q \,\cong \sqcup_2\, q\colon \Bd T(Y,X)\cong \sqcup_2 \, Y \darrow{d} \Bd T(S^3, S^1 \simtimes S^3) \cong \sqcup_2 \, S^3,\]
whose restriction over each copy of $S^3$ is ribbon fillable. Since there is a homeomorphism 
\[(S^1 \simtimes S^3)\setminus \Int(T(S^3, S^1 \simtimes S^3))\cong S^3 \times D^1 \cong S^4 \setminus \Int(\sqcup_2\, D^4),\]
we can apply Theorem \ref{thm:pz} with $n=2$ to extend $\widetilde q$ and hence $\widehat q$ to a simple $d$-fold branched covering
\[p\colon X=T(Y,X)\cup_{\Bd} X' \darrow{d} S^1 \simtimes S^3 \cong T(S^3, S^1 \simtimes S^3)\cup_{\Bd} (S^3 \times D^1)\]
satisfying the condition (a) or (b). Here, the gluing map of the last boundary union is the disjoint union of the identity map and a reflection of $S^3$.
\hfill \qed

\section{Final examples}

\begin{example}
For any $g \in \N$, denote by $\Sigma_g$ (resp. $N_g$) the closed connected orientable (resp. non-orientable) surface of genus $g$. Let $g_1,g_2 \in \N$ and consider the cartesian product $\Sigma_{g_1}\times N_{g_2}$. Represent $\Pd (w_1(N_{g_2}))$ via a simple loop $\gamma \subset N_{g_2}$, which is $2$-sided if and only if $g_2$ is even. We then have that \[[\Sigma_{g_1}\times \gamma]=\Pd(w_1(\Sigma_{g_1}\times N_{g_2})).\] This implies that 
\begin{equation*}
\Pd(w_1(\Sigma_{g_1}\times N_{g_2})^2)=[\Sigma \times \gamma]^2=
    \begin{cases}
      0 & \text{ if $g_2$ is even}\\
      [\Sigma_{g_1}\times \{*\}] & \text{ if $g_2$ is odd}
    \end{cases}\,
\end{equation*}

    while
\[\Pd(w_1(\Sigma_{g_1}\times N_{g_2})^4)=0.\]

For every even $d\geq 4$, Theorem \ref{ThmA} then implies the existence of a simple $d$-fold covering 
\[\Sigma_{g_1}\times N_{g_2}\darrow{d} \RP^4\]
branched over a surface which is nodal if $d=4$ and locally flat if $d \geq 6$. Moreover, we can apply Theorem $\ref{ThmB}$ to deduce that, given $d \geq 4$, there exists a simple $d$-fold branched covering
\[\Sigma_{g_1}\times N_{g_2}\darrow{d} S^1 \simtimes S^3\]
if and only if $g_2$ is even, where the branch set is a nodal surface if $d=4$ and non-singular if $d \geq 5$.
\end{example}

\begin{example}
Let $g_1, g_2 \geq 1$ and consider the product $N_{g_1}\times N_{g_2}$. We have that
\[\Pd(w_1(N_{g_1}\times N_{g_2}))=[N_{g_1}\times \gamma_2]+[\gamma_1\times N_{g_2}]\]
where $\gamma_i \subset N_{g_i}$ is a simple loop dual to $w_1(N_{g_i})$ for $i=1,2$. In particular, $\gamma_i$ is $2$-sided if and only if $g_i$ is even. It follows that
\begin{align*}
{\Pd} & (w_1(N_{g_1}\times N_{g_2})^2)=[N_{g_1}\times \gamma_2]^2+[\gamma_1\times N_{g_2}]^2=\\
   &\begin{cases}
      0 & \text{ if $g_1$ and $g_2$ are even}\\
      [\{*\}\times N_{g_2}] & \text{ if $g_1$ is odd and $g_2$ is even}\\
      [N_{g_1}\times \{*\}] & \text{ if $g_1$ is even and $g_2$ is odd}\\
      [N_{g_1}\times \{*\}]+[\{*\}\times N_{g_2}] & \text{ if $g_1$ and $g_2$ are odd}
    \end{cases}
\end{align*}
and hence
\begin{equation*}
    \Pd(w_1(N_{g_1}\times N_{g_2})^4)=0.
\end{equation*}

As a consequence of Theorem \ref{ThmA}, for every even $d\geq 4$ there exists a simple $d$-fold covering
\[N_{g_1}\times N_{g_2}\darrow{d}\RP^4\]
branched over a surface which is nodal if $d=4$ and non-singular if $g \geq 6$.

Moreover, Theorem \ref{ThmB} implies that, given $d \geq 4$, there exists a simple $d$-fold branched covering
\[N_{g_1}\times N_{g_2} \darrow{d} S^1 \simtimes S^3\]
if and only if both $g_1$ and $g_2$ are even, where the branch set is again a nodal surface if $d=4$ and it is non-singular if $d \geq 5$.
\end{example}
We conclude with the following observation about exotic projective $4$-spaces.
\begin{example}
The conditions in Theorem \ref{ThmA} and \ref{ThmB} implying the existence of simple branched coverings over $\RP^4$ and $S^1 \simtimes S^3$ are purely topological. In particular, they do not depend on the PL structure of the covering $4$-manifold $X$. As a consequence, we have the following. For every odd $d \geq 5$ and every $4$-manifold $\RP^4_{\!\!\textup{ex}}$ homeomorphic but not PL diffeomorphic to $\RP^4$, there exists a simple $d$-fold covering $\RP^4_{\!\!\textup{ex}} \rightarrow \RP^4$ branched over a locally flat PL surface. For example, we can take $\RP^4_{\!\!\textup{ex}}$ to be the exotic $\RP^4$ constructed in \cite{[CappellShaneson]} or in \cite{[FintushelStern1]}.
\end{example}


\begin{thebibliography}{99}

\bibitem{A1920} J.\,W.~Alexander, \textsl{Note on Riemann spaces}, Bull. Amer. Math. Soc. {\bf 26} (1920), no. 8, 370--372.


\bibitem{BE} I.~Berstein and A.\,L.~Edmonds, \textsl{On the construction of branched coverings\break of low-dimensional manifolds}, Trans. Amer. Math. Soc. {\bf 247} (1979), 87--124.

\bibitem{BP} I.~Bobtcheva and R.~Piergallini, \textsl{On 4-dimensional 2-handlebodies and 3-mani\-folds}, J. Knot Theory Ramifications {\bf 21} (2012), no.~12, 1250110, 230 pp.

\bibitem{[CappellShaneson]} S.~E. Cappell and J.~L. Shaneson, \textsl{Some new four-manifolds}, Ann. of Math. (2) {\bf 104} (1976), no.~1, 61--72.

\bibitem{E} A.\,L.~Edmonds, \textsl{Extending a branched covering over a handle}, Pacific J. Math. {\bf 79} (1978), no.~2, 363--369.

\bibitem{[FintushelStern1]}  R. Fintushel and R. J. Stern, \textsl{An exotic free involution on $S^4$}, Ann. of Math. 113 (1981), 357.

\bibitem{Hi} H.~M. Hilden, \textsl{Every closed orientable $3$-manifold is a $3$-fold branched covering space of $S\sp{3}$}, Bull. Amer. Math. Soc. {\bf 80} (1974), 1243--1244.

\bibitem{HL} H.\,M.~Hilden and R.\,D.~Little, \textsl{Cobordism of branched covering spaces}, Pacific J. Math., {\bf 87} (1980), 335--345.

\bibitem{Hirsch} U. Hirsch, \textsl{\"Uber offene Abbildungen auf die 3-Sph\"are}, Math Z. {\bf 140} (1974), 203--230.

\bibitem{M} J.\,M.~Montesinos-Amilibia, \textsl{$4$-manifolds, $3$-fold covering spaces and ribbons}, Trans. Amer. Math. Soc. {\bf 245} (1978), 453--467.

\bibitem{M1} J.\,M. Montesinos, \textsl{Lectures on 3-fold simple coverings and 3-manifolds}, Contemp. Math. {\bf 44} (1985), 157-177.

\bibitem{M2} J.\,M. Montesinos, \textsl{A note on moves and irregular coverings of $S^4$}, Contemporary Mathematics {\bf 44} (1985), 345–349.

\bibitem{P} R.~Piergallini, \textsl{Four-manifolds as 4–fold branched covers of $S^4$}, Topology {\bf 34} (1995) 497–508.

\bibitem{PI} R.~Piergallini and M.~Iori, \textsl{4-manifolds as covers of the 4-sphere branched over non-singular surfaces}, Geometry \& Topology {\bf 6} (2002), 393--401.

\bibitem{PZ1} R.~Piergallini and D.~Zuddas, \textsl{On branched covering representation of 4-mani\-folds}, J. Lond. Math. Soc. (2) {\bf 100} (2019), no.~1, 1--16.

\bibitem{PZ2} R.~Piergallini and D.~Zuddas, \textsl{Branched coverings of $\CP^2$ and other basic\break 4-manifolds}, Bull. Lond. Math. Soc. {\bf 53} (2021), no.~3, 825--842.

\bibitem{T} R.~Thom, \textsl{Quelques propriétés globales des variétés diffêrentiables}, Comment. Math. Helv. {\bf 28} (1954), 17--86.

\end{thebibliography}
\end{document}